\documentclass[11pt]{amsart}
\usepackage{amsmath}
\usepackage[dvips]{graphicx}
\usepackage{amsfonts}
\usepackage{amsopn}
\usepackage{amsthm} 
\usepackage{amssymb}
\usepackage{color}
\usepackage{mathrsfs}
\usepackage{hyperref}

\usepackage{latexsym}
\usepackage{amsgen}
\usepackage{comment}

\usepackage[title]{appendix}
\usepackage{algorithm}
\usepackage{algpseudocode}

\DeclareMathOperator*{\argmax}{arg\,max}
\DeclareMathOperator*{\argmin}{arg\,min}
\newtheorem{theorem}{Theorem}[section]

\newtheorem{problem}[theorem]{Problem}

\numberwithin{equation}{section}

\def\Om{\Omega}

\newcommand{\be}{\begin{equation}}
\newcommand{\ee}{\end{equation}}

\newcommand{\ba}{\begin{array}}
	\newcommand{\ea}{\end{array}}

\usepackage{amssymb}

\numberwithin{equation}{section}

\title[Inverse problem in mean field games from partial boundary measurement]{A numerical algorithm for inverse problem from partial boundary measurement arising from mean field game problem}

\author{Yat Tin Chow}
\thanks{Department of Mathematics, University of California, Riverside (yattinc@ucr.edu)}

\author{Samy Wu Fung}
\thanks{Department of Applied Mathematics and Statistics, Colorado School of Mines (swufung@mines.edu)}
\author{Siting Liu}
\thanks{Department of Mathematics, University of California, Los Angeles (siting6@math.ucla.edu, lnurbek@math.ucla.edu, sjo@math.ucla.edu)}
\author{Levon Nurbekyan}
\author{Stanley Osher}
\thanks{S. W. Fung, S. Liu, L. Nurbekyan and S. Osher thank the funding from AFOSR MURI FA9550-18-1-0502 and ONR grants: N00014-18-1-2527, N00014-20-1-2093, and N00014-20-1-2787.}
\thanks{Authors are listed in alphabetical order.}

\begin{document} 
	\maketitle
	
	\begin{abstract}
    In this work, we consider a novel inverse problem in mean-field games (MFG). We aim to recover the MFG model parameters that govern the underlying interactions among the population based on a limited set of noisy partial observations of the population dynamics under the limited aperture. 
    Due to its severe ill-posedness, obtaining a good quality  reconstruction is very difficult.  
    Nonetheless, it is vital to recover the model parameters stably and efficiently in order to uncover the underlying causes for population dynamics for practical needs.
		
    Our work focuses on the simultaneous recovery of running cost and interaction energy in the MFG equations from a \emph{finite number of boundary measurements} of population profile and boundary movement. 
    To achieve this goal, we formalize the inverse problem as a constrained optimization problem of a least squares residual functional under suitable norms. 
    We then develop a fast and robust operator splitting algorithm to solve the optimization using techniques including harmonic extensions, three-operator splitting scheme, and primal-dual hybrid gradient method.
    Numerical experiments illustrate the effectiveness and robustness of the algorithm.
	\end{abstract}

\section{Introduction}	
The basis for the MFG framework is the concept of Nash equilibrium, where agents cannot unilaterally improve their objectives. Under suitable regularity assumptions, a common MFG model reduces to the following system of partial differential equations (PDE):
	\begin{eqnarray}
	    \hspace{10pt}
		\begin{cases}
			&-\partial_t \phi(x,t)  - \nu \Delta \phi (x,t)+ H (x, \nabla_x \phi (x,t)) = F (x,\rho(\cdot,t)),~\text{in}~\Omega'\times (0,1),\\
			&\partial_t \rho(x,t)-\nu  \Delta \rho (x,t)- \nabla_x \cdot ( \rho(x,t) \, \nabla_p H(x, \nabla_x \phi (x,t)))  =0,~\text{in}~\Omega'\times (0,1),\\
			&\rho(x,0) = {\rho}_0(x),~\phi(x,1) = g(x),~\text{in}~ \Omega'.
		\end{cases}
		\label{eq:mfg_gen}
	\end{eqnarray}
	Here, $\rho(\cdot,t),~t\in [0,1]$ represents the population distribution over the state space $\Omega'$ at time $t$ satisfying a Fokker-Planck equation, and $\phi(x,t)$ is the value function of each player that satisfies a Hamilton-Jacobi equation and governs optimal actions of players. 
	The Hamiltonian, $H$, is the Legendre transform of the Lagrangian, $L$, representing the running cost for each agent. Furthermore, $F$ represents an interaction between the agents and the population. Typical choices for $H,L,F$ in crowd motion applications are
	\begin{align}
	H(x,p) = \frac{1}{2} \kappa(x) |p|^2\, , \;	L(x,v) = \frac{1}{2 \kappa(x)} |v|^2 \, , \;F(x,\rho) = \int_{\Omega'}  K(x,y) \rho(y) dy \,.
	\label{eq:simplification}
	\end{align}

Introduced in \cite{LasryLions2007,LasryLions06a,LasryLions06b} and \cite{HCM06,HCM07},
MFG is an actively growing field
significantly advancing the understanding of social cooperation and economics \cite{moll14,cardialiaguet2018,gomes2015economic}, biological systems \cite{stella2018mean}, election dynamics \cite{yang2018learning}, population games \cite{lauriere2020convergence}, robotic control \cite{liu2018mean}, machine learning \cite{ruthotto2020machine, lin2021alternating}, dynamics of multiple populations \cite{cirant2015multi}. Recently, they are utilized to understand pandemic modeling and control such as COVID-19 \cite{lee2021controlling}.

With the significant descriptive power of MFGs, it is vital to consider inverse problems arising in MFGs. We aim to reconstruct MFG parameters for a class of nonlocal problems, including the geometry of the underlying space and the interactions between large crowds, based on partial population observations. More specifically, we are interested in the following problem.
\begin{problem}\label{pbm:main}
Given a part of the solution to an MFG system \eqref{eq:mfg_gen}, \eqref{eq:simplification}, for instance,
\begin{equation*}
\Big(\rho\left(x,s\right), -\rho \left(x,s\right)\nabla_p H(x, \nabla_x \phi (x,t))\Big)|_{\partial \Omega \times (0,T)},
\end{equation*}
for finitely many examples of $\rho_0$ and terminal cost $g$, can we numerically recover the speed field $\kappa(x)$ and the interaction kernel $K(x,y)$ from observations?
\end{problem}

Such a model-recovery algorithm can help understand the underlying population dynamics in numerous problems, such as migration flow or contagious rate of COVID-19. We further envision applications to include rescue and exploration team management, policymaking, diplomacy, election modeling, catastrophe management, and evacuation planning.

{Note that $m(x,s)=-\rho\left(x,s\right)\nabla_p H(x, \nabla_x \phi (x,s))$ represents the flux of the agents through the state $x$ at time $s$ as a result of optimal actions. The interpretation of the flux is straightforward for crowd-dynamics models and can be measured by counting people crossing checkpoints or parts of the border. For such models, the value function $\phi(x,s)$ could represent the travel cost for a traveller who is at location $x$ at time $s$. Hence, one could also consider an inverse problem where one observes the value function, instead of the flux, by looking at travel companies' expenses or consumer ticket prices (discounted for the companies' profit margins).

For economic and finance models~\cite{gomes_book16,moll14,cardialiaguet2018} the state variables typically represent asset (wealth, income, inventory) levels instead of a physical location. Hence, the value function represents maximal utility for agents with a given asset level, and the flux represents the total amount of transactions performed by them. Interestingly, in economic models one often has \textit{implicit} mean-field interactions reflected in \textit{market-clearing} type conditions instead of an explicit interaction functional $F(x,\rho)$. Hence, a related inverse problem is to find an appropriate market-clearing condition or tune its parameters for a given economy. This manuscript addresses explicit models with flux data leaving the implicit ones with other data types for future work.

\subsection{Our contributions}

We propose a new MFG inverse problem with non-invasive partial boundary measurements. Based on insights from \cite{nursaude18,liu2020computational,liu2020splitting} we postulate a feature expansion representation for the interaction kernel $K$ and formulate the forward problem \eqref{eq:mfg_gen} as a convex-concave saddle point problem. Furthermore, we design a three-operator splitting scheme~\cite{davis2017three} for the resulting inverse problem with a saddle-point constraint. The algorithm reduces to a forward-backward splitting for the parameter updates, a primal-dual hybrid gradient for the forward saddle point problem update, and a proximal-point algorithm for the adjoint problem update. Intriguingly, our algorithm applies to inverse problems whose forward problem has a saddle point structure beyond MFG.

\subsection{Related work}

Despite of the large body of work on theory, numerical methods, and applications~\cite{mfgCIME}, inverse problems arisen from MFG is still quite an unexplored terrain. To the best of our knowledge, only~\cite{ding2020mean, kachroo2015inverse,Burger_2020} study such problems. The work in \cite{ding2020mean} is the closest to our objective but considers the case with a full space-time measurement of data in the sampling domain. However, most inverse problems in practice only have partial boundary measurements available, either obtained via non-invasive measurement methods or because of the limited access to the sampling domain. Compared with the case with full space-time measurement in the domain, inverse problems with only partial boundary measurements are generally known to be more severely ill-posed.  In this work, we focus on the recovery problem with only boundary measurements coming from several measurement events.

The rest of the paper is organized as follows. In Section \ref{sec:inv_saddle_algo}, we introduce an abstract inverse problem with a saddle point constraint and a generic algorithm to solve it. In Section \ref{sec:inv_mfg}, we present the inverse MFG formulation. Next, in Section \ref{sec:algo_implementation} we discuss the implementation of the algorithm in Section \ref{sec:inv_saddle_algo} for the inverse MFG in Section \ref{sec:inv_mfg}. Section \ref{sec:num_examples} contains three numerical examples to demonstrate the robustness and effectiveness of our algorithm. Finally, Section \ref{sec:conclusion} contains a discussion and concluding remarks.

\section{An inverse problem with a saddle point forward model}\label{sec:inv_saddle_algo}

In this section, we formulate an abstract inverse problem with a saddle point forward model. We discuss suitable Karush–Kuhn–Tucker (KKT) conditions and a generic algorithm to solve such inverse problems.

\subsection{A forward saddle point problem}
Consider a saddle point problem
   \begin{equation}\label{eq:infsup_F}
       \min_{x \in \mathcal{X}} \max_{y \in \mathcal{Y}} F(u,x,y,c),
   \end{equation}
	where $F : \mathcal{U} \times X \times Y \times \mathcal{D} \rightarrow \mathbb{R} \bigcup \{\pm \infty\}$ a smooth functional such that $(x,y) \mapsto F(u,x,y,c)$ is strongly convex-concave. 
	Here, $x$ is the primal variable, and $y$ is the dual variable in the forward problem. Commonly, $y$ is either used to handle constraints in the forward problem or linearize nonlinear components via some splitting scheme. The variable $c$ represents model parameters associated with the functional $F$, while $u$ represents boundary and initial-terminal conditions. Given model parameters $c$, we define a boundary measurement map $\Lambda_{c}$ as follows:
	\begin{eqnarray*}
		\Lambda_{c} : \mathcal{U} &\rightarrow& \Pi_{B} X \times \Pi_{B} Y  \\ 
		u &\mapsto& \left(\Pi_{B,x} \left(x\right), \Pi_{B,y} \left(y\right)\right) \; \text{where} \; (x,y) \in \argmin_{x} \argmax_{y}  F(u,x,y,c),
	\end{eqnarray*}
	where $\Pi_{B,x}, \Pi_{B,y}$ denote a projection operator that represent the partial boundary measurements of $x,y$ available. 
	We note that $u$ corresponds to boundary conditions of the forward problem, whereas $B$ is the subset of the domain where the partial measurements are collected.

	\subsection{The inverse problem and a generic algorithm}

	Assume that 
	\begin{eqnarray*}
		( \tilde{r}_{B,i}, \tilde{s}_{B,i})  \approx \Lambda_{c} (u_i) = \left(  [\Pi_{B,x} \left(x\right)]\left(u_i\right),   [\Pi_{B,y} \left(y\right)]\left(u_i\right) \right)
	\end{eqnarray*}
	are noisy measurements for a given $\{u_i\}_{i=1}^N \in \mathcal{U}$. Our goal is to recover $c \in \mathcal{D}$. We formulate this problem as a constrained optimization problem
	\begin{equation}
	   \label{prob:inv_saddle}
	   \begin{aligned}
		\inf_{\{x_i,y_i\}_{i=1}^N,c} \bigg \{ \sum_{i=1}^N \frac{1}{2} \| \Pi_{B,x} (x_i) - \tilde{r}_{B,i}  \|^2 +  \sum_{i=1}^N  \frac{1}{2} \| \Pi_{B,y} (y_i) - \tilde{s}_{B,i}  \|^2 + R(c) \, : \\  (x_i,y_i) \in \argmin_{x} \argmax_{y}  F(u_i,x,y,c) \bigg\},
	\end{aligned}
	\end{equation}
	where $R$ is a suitable regularizer and $\|\cdot\|$ are suitable choices of (semi)-norms. Introducing Lagrange multipliers (dual variables) $(\lambda_{x_i}, \lambda_{y_i} )$, \eqref{prob:inv_saddle} reduces to
	\begin{equation}\label{prob:inv_saddle_saddle_form}
	\begin{aligned}
		& &\inf_{\{x_i,y_i\}_{i=1}^N,c}  \sup_{\{ \lambda_{x_i}, \lambda_{y_i} \}_{i=1}^N }  \bigg \{  \sum_{i=1}^N \frac{1}{2} \| \Pi_{B,x} (x_i) - \tilde{r}_{B,i}  \|^2 + \sum_{i=1}^N  \frac{1}{2} \| \Pi_{B,y} (y_i) - \tilde{s}_{B,i}  \|^2 + R(c) \\
		&  & \quad \quad  \quad \quad   \quad \quad \quad \quad \quad  \quad + \sum_{i=1}^N  \langle \partial_{x_i} F(u_i,x_i,y_i,c) , \lambda_{x_i}  \rangle -  \sum_{i=1}^N   \langle \partial_{y_i} F(u_i,x_i,y_i,c) , \lambda_{y_i}  \rangle  \bigg\}.
	\end{aligned}
	\end{equation}
	Thus, the KKT condition for \eqref{prob:inv_saddle}, \eqref{prob:inv_saddle_saddle_form} are as follows:
		\begin{equation}
		\label{eq:KKT_conditions}
		\begin{aligned}
		\begin{cases}
			\Pi_{B,x}^* [ \Pi_{B,x} (x_i) - \tilde{r}_{B,i}  ] + \partial^2_{x_i,x_i} F(u_i,x_i,y_i, c)  \lambda_{x_i}  &= 0, \\
			\Pi_{B,y}^* [ \Pi_{B,x} (y_i) - \tilde{s}_{B,i}  ] - \partial^2_{y_i,y_i} F(u_i,x_i,y_i, c)  \lambda_{y_i}  &= 0,\\
			\partial_c R(c) + \sum\limits_{i=1}^N  \langle \partial_c \partial_{x_i} F(u_i,x_i,y_i,c) , \lambda_{x_i}  \rangle - \sum\limits_{i=1}^N  \langle \partial_c \partial_{y_i} F(u_i,x_i,y_i,c) , \lambda_{y_i}  \rangle   &= 0, \\
			\partial_{x_i} F(u_i,x_i,y_i,c) &= 0,\\
			- \partial_{y_i} F(u_i,x_i,y_i,c) &= 0,
		\end{cases}
		\end{aligned}
		\end{equation}
	for $i=1,\ldots,N$. Here, $\Pi_{B,x}^*, \Pi_{B,y}^*$ are the adjoints of $\Pi_{B,x}, \Pi_{B,y}$, respectively. 

Finally, we formulate these KKT conditions as an inclusion problem
	\begin{equation*}
	\begin{aligned}
	0 \in A\left( c ,  (x, y), (\lambda_x,\lambda_y)  \right)  + C\left( c, (x, y), (\lambda_x,\lambda_y)  \right),
	\end{aligned}
	\end{equation*}
	where
	\begin{equation*}
	\begin{aligned}
	A\left( c , (x, y), (\lambda_x,\lambda_y)\right) = 
	\begin{pmatrix}
	\partial_c R(c) \\
	(0,0)\\
	(0,0) \\
	\end{pmatrix},
	\end{aligned}
	\end{equation*}
	and
	\begin{eqnarray*}
		&  &C\left( c , (x, y), (\lambda_x,\lambda_y)\right)  \\
		&=&
		\begin{pmatrix}
			\sum_{i=1}^N  \langle \partial_c \partial_{x_i} F(u_i,x_i,y_i,c) , \lambda_{x_i}  \rangle - \sum_{i=1}^N  \langle \partial_c \partial_{y_i} F(u_i,x_i,y_i,c) , \lambda_{y_i}  \rangle  \\
			\big(  \partial_{x_i} F(u_i,x_i,y_i,c),  - \partial_{y_i} F(u_i,x_i,y_i,c) \big) \\
			\big(  \partial^2_{x_i,x_i} F(u_i,x_i,y_i, c)  \lambda_{x_i}  + \Pi_{B,x}^* [ \Pi_{B,x} (x_i) - \tilde{r}_{B,i}  ]  , \\  - \partial^2_{y_i,y_i} F(u_i,x_i,y_i, c)  \lambda_{y_i}  + \Pi_{B,y}^* [ \Pi_{B,x} (y_i) - \tilde{s}_{B,i}  ] \big)  \\
		\end{pmatrix}.
	\end{eqnarray*}
Note that $A$ is monotone but $C$ is not known to be monotone in general. 
	

\subsection{A generic algorithm}\label{alg:inv_saddlepoint}
	
Here, we outline an iterative algorithm for solving \eqref{prob:inv_saddle_saddle_form}. At $(n+1)$-th iteration, we first update the adjoint variables $\{(\lambda_{x_i},\lambda_{y_i})\}_{i=1}^N$ using the Chambolle-Pock method~\cite{champock11}; then we update $c$ for the inverse problem by taking a proximal gradient step; next we use the Chambolle-Pock method again to compute forward problems $\{(x_i,y_i)\}_{i=1}^N$. Summarizing, a high level description of the $(n+1)$-th iteration is as follows: 
\begin{eqnarray*}
\begin{cases}
\begin{cases}
		\lambda_{y_i}^{n+1}  &= [ 1 - \alpha_{\lambda_{y_i}}  \partial^2_{y_i,y_i}  F(u_i,x_i^n,y_i^n, c^n) ]^{-1}  \left(  \lambda_{y_i}^n  -  \alpha_{\lambda_{y_i}}   \Pi_{B,y}^* [ \Pi_{B,y} (y_i^n) - \tilde{s}_{B,i}  ]  \right) \\
		\lambda_{x_i}^{n+1,\text{temp}}  &=  [  1+  \alpha_{\lambda_{x_i}} \partial^2_{x_i,x_i}  F(u_i,x_i^n,y_i^n, c^n) ]^{-1}  \left( \lambda_{x_i}^n  - \alpha_{\lambda_{x_i}}  \Pi_{B,x}^* [ \Pi_{B,x} (x_i^n) - \tilde{r}_{B,i}  ] \right) \\
	        \lambda_{x_i}^{n+1} & = 2 \lambda_{x_i}^{n+1,\text{temp}} - \lambda_{x_i}^{n,\text{temp}}
\end{cases} 
 \\
\begin{cases}
		c^{n+1}  &= ( I + \alpha_{c} \partial_c R )^{-1} \Big[ c^{n} - \alpha_{c}  \sum_{i=1}^N  \langle \partial_c \partial_{x_i}  F(u_i,x_i^n,y_i^n,c^n) , \lambda^{n+1}_{x_i}  \rangle \\& \quad  + \alpha_c \sum_{i=1}^N  \langle \partial_c \partial_{y_i}  F(u_i,x_i^n,y_i^n,c^n) , \lambda^{n+1}_{y_i}  \rangle \Big]
\end{cases} 
\\
\begin{cases}
		x^{n+1}_i&= [ 1 + \alpha_{x_i} \partial_{x_i}  F(u_i, \cdot ,y_i^n,c^{n+1}) ]^{-1} ( x^n_i ) \\
    		y^{n+1,\text{temp}}_i &=  [ 1 - \alpha_{y_i}  \partial_{y_i}  F(u_i,x_i^{n+1},\cdot ,c^{n+1}) ]^{-1} (y^n_i) \\
        y^{n+1} & = 2 y^{n+1,\text{temp}}_i - y^{n,\text{temp}}_i,
\end{cases}

\end{cases}
\end{eqnarray*}
where  $ (\alpha_{\lambda_{x_i}}, \alpha_{\lambda_{y_i}},\alpha_c,  \alpha_{x_i},\alpha_{y_i} )$ are the corresponding time steps.

In what follows, we specify the MFG inverse problem and the implementation of the algorithm above for it.

\section{An inverse MFG problem}\label{sec:inv_mfg}

Here, we explain the saddle point problem formulation of nonlocal MFG~\cite{nursaude18,liu2020computational,liu2020splitting} and formulate the inverse MFG problem of our interest.

\subsection{Saddle point formulation of MFG via feature-space expansions}

Consider the following MFG system with nonlocal couplings:
\begin{equation} \label{eq:nonlocal_mfg}
\begin{cases}
-\phi_t(x,t)-\nu \Delta \phi(x,t)  +  \frac{\kappa(x)}{2}  \|\nabla \phi(x,t)\|^2 = \int_{\Om'} K(x,y) \rho(y,t) dy  &\text{ in } \Omega' \times (0,1),\\
\rho_t(x,t)  -\nu \Delta \rho(x,t)-\nabla \cdot (\kappa(x)\rho(x,t) \nabla \phi(x,t))=0 & \text{ in } \Omega'\times (0,1),\\
(\kappa(x)\rho(x,t) \nabla \phi(x,t)) \cdot n = 0 &  \text{ on } \partial \Omega'\times (0,1),\\
\rho(x,0)=\rho_0(x),~\phi(x,1)=g(x) &  \text{ in } \Omega'.\\
\end{cases}
\end{equation}

We assume that $K$ is positive definite and translation invariant, which yields that the mean-field interaction satisfies the Lasry-Lions monotonicity condition~\cite{LasryLions2007} and agents are crowd averse. Moreover, \eqref{eq:nonlocal_mfg} admits a saddle point formulation
\begin{equation}\label{eq:nonlocal_infsup_direct}
	\begin{split}
		& \inf_{\substack{\phi}} \sup_{\rho,m}   \Bigg\{-\int_{\Omega' }\phi(x,0)\rho_0(x) dx - \int_{\Omega' }\int_0^1 \left(\rho \phi_t  
		+\nu \rho \Delta \phi + m \cdot \nabla \phi \right) dt dx\\
		&-\int_{\Omega'}\int_0^1 \frac{1}{2\kappa(x)}\frac{\|m\|^2}{2\rho}  dt dx-\frac{1}{2} \int_{\Omega'\times \Omega'} K(x,y)\rho(x,t)\rho(y,t)dxdy  - \chi_{\rho \geq 0} + \chi_{\phi(x,1) =g(x)} \Bigg\}\\
	\end{split}
\end{equation}
Here, $\chi_{Z}(z)$ is the indicator function over the set $Z$ defined by
\begin{align*}
	\chi_{Z}(z)=\begin{cases}
		0, \quad \text{if} \; z\in Z\\
		\infty, \quad \text{otherwise}.
	\end{cases}
\end{align*}
Modeling the interaction term
\begin{equation*}
	\frac{1}{2} \int_{\Omega'\times \Omega'} K(x,y)\rho(x,t)\rho(y,t)dxdy
\end{equation*}
directly is costly for both forward model and the inverse problem. Moreover, based on the works from~\cite{nursaude18,liu2020computational,liu2020splitting,agrawal2022random}, we model and approximate this term using feature-space expansions. More specifically, based on Bochner's theorem \cite{rong1998bochner}, we postulate that
\begin{equation*}
\begin{split}
K(x,y)=&\sum_{k=1}^r \mu_k^2 \cos(\omega_k\cdot (x-y))\\
=&\sum_{k=1}^r \left(\mu_k^2 \cos (\omega_k \cdot x) \cos (\omega_k \cdot y)+\mu_k^2 \sin (\omega_k \cdot x) \sin (\omega_k \cdot y) \right)
\end{split}
\end{equation*}
for some $\{\omega_k\} \subset \mathbb{R}^d$, and $\{\mu_k^2\} \subset \mathbb{R}_+$.
Denoting by
\begin{equation*}
	\begin{split}
	\mu =&\big(\mu_{1,1},\mu_{1,2}, \mu_{2,1}, \cdots, \mu_{r,1}, \mu_{r,2}\big) \quad
	\omega = \big(\omega_{1,1},\omega_{1,2}, \omega_{2,1}, \cdots, \omega_{r,1}, \omega_{r,2}\big)\\
	C_{\text{odd}=\text{even}} =& \bigg\{ (x_{1,1},x_{1,2}, x_{2,1}, \cdots, x_{r,1}, x_{r,2}): x_{i,1} = x_{i,2} \bigg\}	\\
	\zeta(x;\mu,\omega) =& \bigg(\mu_{1,1} \cos(\omega_{1,1} \cdot x), \mu_{1,2} \sin(\omega_{1,2} \cdot x),\cdots,\mu_{r,1} \cos(\omega_{r,1} \cdot x),\mu_{r,2} \sin(\omega_{r,2} \cdot x)\bigg) 
	\end{split}
\end{equation*}
we obtain
\begin{equation*}
	K(x,y)=\zeta(x;\mu,\omega)\cdot \zeta(y;\mu,\omega), \quad \mu,\omega \in C_{\text{odd}=\text{even}}.
\end{equation*}
Using this representation, we obtain
\begin{equation*}
	\begin{split}
		\frac{1}{2}\int_{\Omega'\times \Omega'} K(x,y)\rho(x,t)\rho(y,t)dt=& \frac{1}{2} \left\| \int_{\Omega'} \zeta(x;\mu,\omega) \rho(x,t) dx \right\|^2 \\
		=&\sup_{a} \left\{a(t) \cdot \int_{\Omega'} \zeta(x;\mu,\omega) \rho(x,t) dx-\frac{1}{2} \int_0^1 \|a(t)\|^2 dt \right\},
	\end{split}
\end{equation*}
where $a(t) = \left(a_{1,1}\left(t\right),a_{1,2}\left(t\right), .., a_{r,1}\left(t\right), a_{r,2}\left(t\right)\right)$ are auxiliary dual variables.
The last equality is a result from~\cite{nursaude18}.
Hence, \eqref{eq:nonlocal_infsup_direct} transforms to
\begin{equation}\label{eq:nonlocal_infsup}
\begin{split}
& \inf_{\substack{\phi, a}} \sup_{\rho,m}   \Bigg\{ \frac{1}{2} \int_0^1 \|a(t)\|^2dt -\int_{\Omega' }\phi(x,0)\rho_0(x) dx - \int_{\Omega' }\int_0^1 \left(\rho \phi_t  
+\nu \rho \Delta \phi + m \cdot \nabla \phi \right) dt dx\\
&-\int_{\Omega'}\int_0^1 \left( \frac{1}{2\kappa(x)}\frac{\|m\|^2}{2\rho} + \rho~a(t)\cdot \zeta(x;\mu,\omega) \right) dt dx- \chi_{\rho \geq 0} + \chi_{\phi(x,1) =g(x)} \Bigg\}\\
& := \inf_{\substack{\phi, a}} \sup_{\rho,m}   \Bigg\{ -\mathcal{L} \bigg(  (\rho_0, g ),  ( \rho, m), (\phi,a)  , (\kappa,\mu)   \bigg)  \Bigg\},\\
\end{split}
\end{equation} 
For more details on representation of nonlocal MFG interactions via a basis and computational methods, see \cite{nursaude18,liu2020splitting,liu2020computational,agrawal2022random}. \textcolor{black}{We also attach an example Algorithm \ref{alg:forward} for calculating the nonlocal mean-field game problem in the appendix.}

\subsection{An inverse mean-field game problem}

Denoting by
\begin{equation*}
\begin{split}
	 & u= (\rho_{0}, g),\quad x= (\rho,m),\quad y= (\phi,a),\quad c = ( \kappa, \mu ),\\
	& F\big( (\rho_{0}, g), (\rho,m),(\phi,a) , (\kappa, \mu ) \big)
	=  \mathcal{L} \big( (\rho_{0}, g), (\rho,m),(\phi,a) , (\kappa, \mu ) \big),   	 
\end{split}
\end{equation*}
we place the MFG forward model in the abstract framework \eqref{eq:infsup_F}. Next,we  assume that $\Omega \subset \Omega'$ and $\kappa(x)$ is known in the domain $\Omega'\backslash \Omega$. We refer to $\Omega$ and $\Omega'$ as sampling and computational domains, respectively. An example is shown in Figure \ref{fig:example_domain}, where the $\Omega'$ is the large square domain, while $\Omega$ is the inner square with its boundary highlighted in red.

	\begin{figure}[h]
	\includegraphics[width=0.950\textwidth,trim=10 10 20 20, clip=flase]{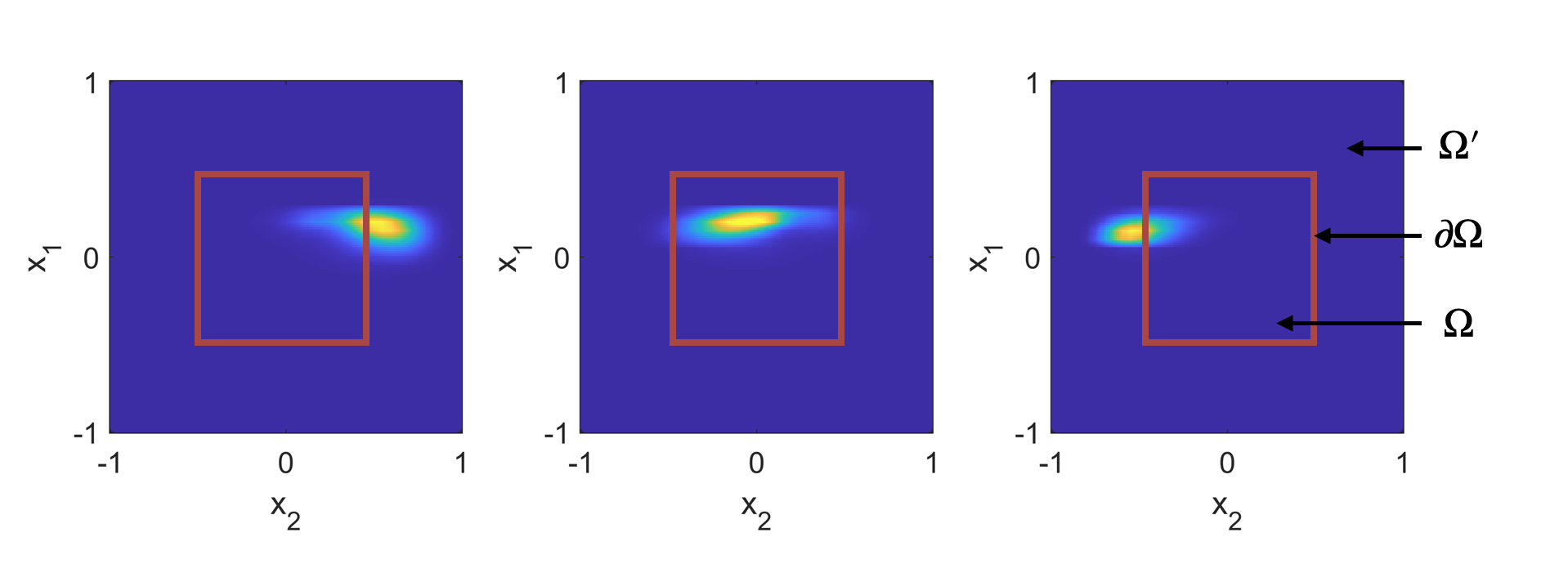}
	\caption{Denote $\rho$ as the solution to the mean-field game system. From left to right, the pictures display the density distribution $\rho$ at time $t = 0.1,  0.5, 0.9$.
	The solid red line represents the boundary of domain $\Omega$. In this mean-filed game, the density travels from the right towards the left, crossing the boundary $\partial \Omega$ twice.}
	\label{fig:example_domain}
	\end{figure}

Next, we take
\begin{align*}
    \Pi_{B,(\rho,m)} \left(\rho,m\right) := \left(  \rho, m \cdot n \right) |_{[0,T] \times \partial \Omega^+ } ,  \quad   \Pi_{B,(a,\phi)} (a,\phi) := \left( 0,0\right),
\end{align*}
for the partial boundary measurement along the boundary $\partial \Omega$. Here, $\partial \Omega^{+}$ means that the normal vector $n$ is pointing outward. Measuring the density and flux through $\partial \Omega$ is reasonable based on physical meaning of the variables. We cannot measure $a$ directly because it is a non-physical auxiliary variable introduced specifically for an efficient representation of nonlocal interactions.


We assume that the ground truth parameters $(\kappa,\mu)$ represent a disturbance of background parameters $(\kappa_0,\mu_0)$. Therefore, given an additional parameter $\varepsilon \geq 0$, we would also like to have a regularization term in the form of $R_{\beta} = \chi_{\beta \geq \varepsilon} (\beta) $.
We also write $R (\kappa,\mu)= R_{1}(\kappa,\mu) + R_{2}(\kappa,\mu)$, where
\begin{align*}
    & R_1(\kappa,\mu) = \tilde{R}_{1}(\kappa) +  \tilde{R}_{2}(\mu) := \gamma_{\kappa} \|\kappa  - \kappa_0\|_{L^1} +  \gamma_{\mu}  \|\mu - \mu_0\|_{L^1},\\
    & R_2(\kappa,\mu):= \chi_{\kappa \geq \varepsilon_1} (\kappa) +   \chi_{\mu^2 \leq \varepsilon_2} (\mu).
\end{align*}

 It is also possible to have other choices of regularization for $(\kappa,\mu)$, such as $TV, H^1$, Wavelet norms.

Now, we can formulate the inverse MFG as follows:
	\begin{equation}
	   \label{prob:inv_mfg}
	   \begin{aligned}
		\inf_{\substack{\{(\rho_i,m_i),(\phi_i,a_i)\}_{i=1}^N,\\\kappa,\mu}} \Bigg \{ \sum_{i=1}^N \frac{1}{2} \| \Pi_{B} (\rho_i,m_i) - \tilde{r}_{B,i}  \|^2 + R (\kappa,\mu) \, : \\  (\rho_i,m_i),(\phi_i,a_i)\in \argmin_{\rho,m} \argmax_{\phi,a}  F((\rho_{0,i}, g_i), (\rho,m),(\phi,a) , (\kappa, \mu )) \Bigg\}.
	\end{aligned}
	\end{equation}


\section{The algorithm}\label{sec:algo_implementation}

We propose an inverse algorithm adapted from the three-operator splitting method~\cite{davis2017three}, which has also been shown to predict Nash equilibria in traffic flows~\cite{heaton2021learn}. We also discuss stabilizing techniques that are essential in practice.

\subsection{The three-operator splitting scheme}
Denoting by $\lambda_{(\rho,m)} := ( \lambda_{\rho}, \lambda_{m} )$ and $\lambda_{(a,\phi)} := ( \lambda_{a}, \lambda_{\phi} )$ and applying the framework in Section \ref{sec:inv_saddle_algo} to \eqref{prob:inv_mfg} we obtain the following inclusion formulation of the inverse MFG problem:
\begin{eqnarray}
\label{eq:monotone_inv_mfg}
0 &\in & A\left( \kappa, \mu\right) + B\left(  \kappa, \mu \right)  + C\left( (\kappa, \mu) ,  ( ( \rho, m) ,(a,\phi)), (\lambda_{(\rho,m)}, \lambda_{(a,\phi)} ) \right) ,
\end{eqnarray}
where 
\begin{eqnarray*}
A\left( \kappa, \mu\right)= \begin{pmatrix}  \partial R_1 ( \kappa, \mu) \\ (0,0) \\ (0,0) \end{pmatrix} \,,\quad 
B\left( \kappa, \mu\right) = \begin{pmatrix} \partial R_2 ( \kappa, \mu) \\ (0,0) \\ (0,0) \end{pmatrix} \,, \\
\end{eqnarray*}
and

\begin{eqnarray*}
&  &C\left( ( \kappa, \mu) ,  ( ( \rho, m) ,(a,\phi)), (\lambda_{(\rho,m)}, \lambda_{(a,\phi)} ) \right)  \\
&=&
\begin{pmatrix}
\sum_{i=1}^N \left \langle  \partial_{( \kappa, \mu)} \partial_{(\rho_i,m_i)}
 \mathcal{L}, \lambda_{(\rho_i,m_i)} \right \rangle  \\
- \sum_{i=1}^N  \left \langle \partial_{( \kappa, \mu)}  \partial_{(a_i,\phi_i)}  \mathcal{L} , \lambda_{(a_i,\phi_i)} \right \rangle   \\
\left(  \partial_{(\rho_i,m_i)}   \mathcal{L},   -  \partial_{(a_i,\phi_i)}  \mathcal{L} \right) \\
 \bigg( \partial^2_{(\rho_i,m_i),(\rho_i,m_i)}
 \mathcal{L}\lambda_{(\rho_i,m_i)}  +  \Pi_{B,(\rho,m)}^* [ \Pi_{B,(\rho,m)} (\rho_i,m_i) - \tilde{r}_{B,i}  ] ,   \\
\quad 
- \partial^2_{(a_i,\phi_i),(a_i,\phi_i)}  \mathcal{L} \lambda_{(a_i,\phi_i)} 
\bigg)  \\

\end{pmatrix}.
\end{eqnarray*}

The three-operator splitting scheme in \cite{davis2017three} applies to optimization problems of the form
\begin{align}
\label{eq:three_op}
    \text{find} \quad z \in \mathcal{H} \quad \text{such that} \quad 0\in Az + Bz + Cz,
\end{align}
where $A,B,C$ are maximal monotone operators defined on a Hilbert space $\mathcal{H}$, and $C$ is cocoercive. Denote by $I_{\mathcal{H}}$ the identity map in $\mathcal{H}$, and $J_S:= (I+S)^{-1}$ the resolvent of a monotone operator $S$. The splitting scheme for solving \eqref{eq:three_op} can be summarized as follows
\begin{equation}
\begin{aligned}
\label{eq:three_op_scheme}
    z^{k+1} &:= (1-\lambda_k)z^k + \lambda_k Tz^k,\\
    T&:= I_{\mathcal{H}} -J_{\gamma B}+J_{\gamma A} \circ (2 J_{\gamma B} -I_{\mathcal{H}} - \gamma C\circ J_{\gamma B} ),
\end{aligned}
\end{equation}
where $\gamma$ is a scalar.  If an operator $S$ is of the sub-differential forms; that is, $S = \partial f_S$ for some functional $f_S$, the resolvent $J_S$ reduces to the proximal map $x\mapsto \argmin\limits_y f_S(y) + \frac{1}{2} \|x-y\|^2$. 

Overall, the algorithm for \eqref{eq:monotone_inv_mfg} follows three components of the generic framework in Section \ref{alg:inv_saddlepoint}, upon some modification. In what follows, we discuss each component separately.

\subsubsection{Update of the adjoint problem}

Firstly, we choose $\| \cdot \|_{(\rho,m)}$ and $\| \cdot \|_{(a, \phi)}$ with
$L^2([0,T], H^{-1/2}(\partial \Omega)) \times  L^2([0,T], H^{1/2}(\partial \Omega)) $ semi-norm, and
$L^2_{t} \times L^2_{x,t}$ norms, respectively.
Here the $H^{\frac{1}{2}}(\partial \Omega)$ and $H^{-\frac{1}{2}}(\partial \Omega)$ semi-norm are taken as follows:
\begin{eqnarray*}
| v |^2_ { H^{\frac{1}{2}} ( \partial \Omega )} &:=&  \min_{v_0 \in H^1 (\Omega) \,,\,  v_0=v \text{ on } \partial \Omega } | v_0 |^2_ { H^{1} ( \Omega ) } \,, \\
| v |^2_ { H^{-\frac{1}{2}} ( \partial \Omega ) } &:=&  \min_{v_0 \in H^1_0 (\Omega) \,, \, \partial_n v_0= \partial_n v \text{ on } \partial \Omega } | v_0 |^2_ { H^{1} ( \Omega ) } \,, 
\end{eqnarray*}
where the right hand side denotes the standard $H^1(\Omega)$ semi-norm.

Assuming appropriate regularity of $(\rho,m)$, we
recall that the operator $\Pi_{B,(\rho,m)}$ is the restriction/trace operator onto the appropriate Sobolev space on the boundary $L^2 \big([0,T]$, $H^{-1/2}(\partial \Omega)\big) \times  L^2\left([0,T], H^{1/2}(\partial \Omega)\right)$
\begin{align*}
\Pi_{B,(\rho,m)} \left(\rho,m\right) := \left(  \rho, m \cdot n \right) |_{[0,T] \times \partial \Omega^+ }. 
\end{align*}
With the aforementioned choice of the semi-norms, we naturally have the (formal) adjoint of $\Pi_{B,(\rho,m)}$,  $\Pi_{B,(\rho,m)}^*$, as the Dirichlet and Neumann harmonic extension operators by definition; that is,

\begin{align*}( \eta_i, \nabla \xi_i ):= \Pi_{B,(\rho,m)}^* [ \Pi_{B,(\rho,m)} (\rho_i^n, m_i^n) - \tilde{r}_{B,i}  ],
\end{align*}
where  $(\eta_i,\xi_i)$ satisfy
\begin{equation}\label{alg:adjoint_1}
\begin{cases}
\begin{cases}
\Delta  \eta_i^{n+1} = 0 \text{ in } \Omega\\
 \eta_i^{n+1}  =  \rho_i^n  - \text{pr}_1 \tilde{r}_{B,i}  \text{ on } \partial \Omega^+
\end{cases} \\
\begin{cases}
\Delta  \eta_i^{n+1} = 0 \text{ in } \Omega' \backslash \Omega\\
 \eta_i^{n+1}  =  \rho_i^n  - \text{pr}_1 \tilde{r}_{B,i}  \text{ on } \partial \Omega^- \\
 \eta_i^{n+1}  = 0  \text{ on } \partial \Omega'
\end{cases}
 
\end{cases}
\begin{cases}
\begin{cases}
\Delta \xi_i^{n+1} = 0 \text{ in } \Omega\\
\partial_n  \xi_i^{n+1}   =   m_i^n \cdot n - \text{pr}_2 \tilde{r}_{B,i}  \text{ on } \partial \Omega^+
\end{cases} \\
\begin{cases}
\Delta \xi_i^{n+1} = 0 \text{ in } \Omega' \backslash \Omega\\
\partial_n  \xi_i^{n+1}   =   m_i^n \cdot n + \text{pr}_2 \tilde{r}_{B,i}  \text{ on } \partial \Omega^- \\
\partial_n  \xi_i^{n+1}   =  0 \text{ on } \partial \Omega'
\end{cases} 
\end{cases}
\end{equation}
Here we use ${\text{pr}_1},{\text{pr}_2}$ to denote the projection from the noisy data.
The harmonic extension is taken at each time $t\in [0,1]$ independently.
In the implementation, we use a standard finite difference scheme to compute the harmonic extension on spatial grids for each time grid point.
Note that if we assume $\kappa= \kappa_0$ to be known outside of domain $\Omega$, the measurements of $m$ and $\nabla \phi$ are equivalent, as $ - \kappa_0(x) \rho (x) \partial_n \phi(x) = m (x) \cdot n$ on $\partial \Omega^+$. 

We remark that the techniques of harmonic extension have been applied to various other problems, e.g. over point clouds and in machine learning \cite{shi2018harmonic}.

It is clear to see that $\lambda_{\phi}$ is redundant and $\lambda_{a} = 0$ whenever $0 \in A + B + C $.  Hence, we can consider only $C\left( (\kappa, \mu),  ( ( \rho, m) ,(a,\phi)), (\lambda_{(\rho,m)}, (0,0) ) \right)$. In this case, we preform a primal-dual hybrid gradient method for updating  $\lambda_{(\rho,m)}$:
\begin{equation}\label{alg:adjoint_2}
\begin{aligned}
\begin{cases}
&\begin{pmatrix} \lambda^{n+1,\text{temp}}_{\rho_i}  \\ \lambda^{n+1,\text{temp}}_{m_i} \end{pmatrix}=\\
&\left [ I + \alpha^n_{\lambda_{(\rho_i,m_i)}} \rho_i^n \begin{pmatrix}  \frac{1}{\kappa^n} \frac{|m^n_i|^2}{(\rho^n_i)^3} \;\;  - \frac{1}{\kappa^n} \frac{m^n_i}{(\rho^n_i)^2} \\  - \frac{1}{\kappa^n} \frac{m_i^n}{(\rho^n_i)^2}    \;\; \frac{1}{\kappa^n \rho_i^n} I  \end{pmatrix} \right]^{-1}
\left( \begin{pmatrix} \lambda^{n}_{\rho_i}  \\ \lambda_{m_i^{n}} \end{pmatrix} - \alpha^n_{\lambda_{(\rho_i,m_i)}} \rho_i^n  \begin{pmatrix} \eta^{n+1}_i  \\ \nabla \xi^{n+1}_i \end{pmatrix} \right)\\
&\begin{pmatrix} \lambda^{n+1}_{\rho_i}  \\ \lambda^{n+1}_{m_i} \end{pmatrix}
= 2
\begin{pmatrix} \lambda^{n+1,\text{temp}}_{\rho_i}  \\
	\lambda^{n+1,\text{temp}}_{m_i} \end{pmatrix}
-\begin{pmatrix} \lambda^{n, \text{temp}}_{\rho_i}  \\ 
\lambda^{n, \text{temp}}_{m_i} \end{pmatrix} 
\end{cases}
\end{aligned}
\end{equation}

\subsubsection{Update of the inverse problem}\label{sec:alg_details:4-1-2}

In this part, we focus on the update for the inverse problem variables $(\kappa,\mu)$. 
\begin{equation}\label{alg:inv_parameters}
\begin{cases}
\begin{cases}
{\kappa}^{n+1, \text{temp}}  &= S_{\alpha_{\kappa}^n  \gamma}\bigg(  2 \kappa^{n} - \tilde{{\kappa}}^{n} - \alpha_{\kappa}^n  
\sum_{i=1}^N   \Lambda_{\kappa}(\kappa^n,m_i^n,\rho_i^n,\lambda_{\rho_i}^{n+1},\lambda_{m_i}^{n+1} )    - \kappa_0  \bigg) + \kappa_0 \\
 \tilde{{\kappa}}^{n+1}  &=  \tilde{{\kappa}}^{n} +  {\kappa}^{n+1, \text{temp}}  - \kappa^{n} \\ 
\kappa^{n+1}  &=  \max \left\{ \varepsilon_1, \tilde{{\kappa}}^{n+1} \right\}  
\end{cases} \\
\begin{cases}
{\mu}^{n+1,\text{temp}}   &=
 S_{\alpha_{\mu}^n  \gamma}\bigg(  2\Pi^*_{\mu}(\mu^n) - \tilde{\mu}^n- \alpha_{\mu}^n  
\sum_{i=1}^N  \Lambda_{\mu}(\lambda_{\rho_i}^{n+1},a_i^n)  
 - \Pi^*_{\mu}(\mu_0)  \bigg) + \Pi^*_{\mu}(\mu_0)  \\
\tilde{\mu}^{n+1}  & = \tilde{\mu}^{n} + {\mu}^{n+1,\text{temp}} - \mu^{n} \\
\mu^{n+1} & = \Pi_{\mu} \left( \min \left\{\varepsilon_2, \tilde{\mu}^{n+1}\right\} \right) 
 
 \end{cases}

\end{cases} 
\end{equation}
where $S_{\alpha}(r)$ is the shrinkage operator given as $S_{\alpha}(r) = \text{sign}( r ) \max\{ \left|  r  \right| - \alpha , 0 \}$, and
\begin{equation*}
\begin{aligned}
 &   \Lambda_{\kappa}(\kappa,m,\rho,\lambda_{\rho},\lambda_m ) =  \int_0^T  \frac{1}{2 (\kappa)^2} \frac{\|m(\cdot, s)\|^2}{ (\rho(\cdot,s))^2} \lambda_{\rho}(\cdot,s) - \frac{1}{(\kappa)^2} \frac{m(\cdot,s)}{\rho(\cdot,s)} \lambda_{m}(\cdot,s) ds,\\
&\Lambda_{\mu}(\lambda_{\rho},a)  =     \int_0^T \int_{\Omega'}\Lambda_{1} (\omega,y)  \lambda_{\rho}(y,s)  a(s)   dy ds,\\
&\Lambda_1 (\omega,x)
=
\begin{pmatrix}
\begin{pmatrix}  \cos(\omega_1 \cdot x) & 0 \\
0 &  \sin(\omega_1 \cdot x)  
\end{pmatrix} 
&  \dots 
& 0 \\
\vdots 
& \ddots 
& 0 \\
0 
& \dots 
& 
\begin{pmatrix}  \cos(\omega_r \cdot x) & 0 \\
0 &  \sin(\omega_r \cdot x)  
\end{pmatrix} 
\end{pmatrix}. \\
\end{aligned}
\end{equation*}
Since we have the $\mu,\omega \in C_{\text{odd}=\text{even}} := \{ (x_{1,1}, x_{1,2}, x_{2,1}, x_{2,2}, ..., x_{r,1}, x_{r,2} ) : x_{i,1} = x_{i,2} \forall i = 1,..., r\}$, we write the projector $\partial \chi_{C_{\text{odd}=\text{even}}}$ (where we identify $C_{\text{odd}=\text{even}}$ with $\mathbb{R}^r$) as
\begin{equation*}
    \begin{aligned}
    &\Pi_{\mu}: \mathbb{R}^{2r} \rightarrow C_{\text{odd}=\text{even}} \cong \mathbb{R}^r\\
    &(x_{1,1}, x_{1,2}, x_{2,1}, x_{2,2}, ..., x_{r,1}, x_{r,2} ) \mapsto (\frac{x_{1,1}+ x_{1,2}}{2}, \frac{x_{2,1}+ x_{2,2}}{2}, ..., \frac{x_{r,1} + x_{r,2}}{2} ) \\
\text{and its adjoint as }\\
&  \Pi_{\mu}^*: \mathbb{R}^r \rightarrow \mathbb{R}^{2r}\\
  & (x_1, x_2, ..., x_r) \mapsto (x_{1}, x_{1}, x_{2}, x_{2}, ..., x_{r}, x_{r} ).\\
    \end{aligned}
\end{equation*}

\subsubsection{Update of the forward problem}

As for the forward problem, we use primal--dual hybrid gradient method (PDHG)~\cite{champock11} to update $ \left(\left(\rho_i,m_i\right),\left(\phi_i,a_i\right)\right)$ for each event $i,$ for $ 1\leq i \leq N$. The iterative updates contains three parts: firstly a proximal gradient descent step for $(\rho_i,m_i)$ with stepsizes $(\alpha_{\rho_i}^n,\alpha_{m_i}^n)$; then a proximal gradient ascent step for $(\phi_i,a_i)$ of stepsizes $(\alpha_{\phi_i}^n,\alpha_{a_i}^n)$; lastly an extrapolating step for $(\phi_i,a_i)$. Note that we make the choice of norm $\|\phi\|_{H^1_{x,t}}^2 = \|\phi_t\|_{L^2_{x,t}}^2 + \|\nabla_x \phi\|_{L^2_{x,t}}^2$ for $\phi$, based on the General-proximal Primal-Dual Hybrid Gradient (G-prox PDHG) method~\cite{jacobs2019solving} that can be interpreted as a preconditioning step for obtaining a mesh-size-free convergence rate for the algorithm. Overall, the computation for the forward model follows the computational method proposed in \cite{liu2020computational,liu2020splitting}.

\begin{equation}\label{alg:update forward_mfg}
    \begin{cases}
( \rho_i^{n+1}  , m_i^{n+1} ) &= 
\text{argmin}_{(\rho,m)} \bigg \{ \mathcal{L} \bigg ( (\rho_{0,i}, g_i ) , ( \rho, m), (a^n_i,\phi^n_i)  , (\kappa^{n+1},\mu^{n+1})  \bigg) \\
 &+ \frac{1}{2 \alpha_{\rho_i}^n} \|\rho_i^{n} - \rho \|_{L^2_{x,t}}^2  + \frac{1}{2 \alpha_{m_i}^n} \|m_i^{n} - m \|_{L^2_{x,t}}^2  \bigg \}
\\
( \phi^{n+1 {, \text{temp}} }_i , a_i^{n+1 {, \text{temp}} } ) &= \text{argmin}_{(a,\phi)}
 \bigg\{ - \mathcal{L} \bigg ( (\rho_{0,i}, g_i ), ( \rho^{n { +1 }}, m^{n { +1 } }), (a,\phi) ,(\kappa^{n+1},\mu^{n+1})  \bigg)  \\
 &+ \frac{1}{2 \alpha_{\phi_i}^n } \| \phi^{n}_i -  \phi \|_{H^1_{x,t}}^2  + \frac{1}{2 \alpha_{a_i}^n } \| a^{n}_i -  a \|_{L^2_{t}}^2
 \bigg\} \\
  ( \phi^{n+1}_i , a_i^{n+1} )  & = 2 ( \phi^{n+1,  \text{temp} }_i , a_i^{n+1, \text{temp} } )  - ( \phi^{n}_i , a_i^{n} ) 
\end{cases} 
\end{equation}

Assembling all three components described above, we arrive at the following algorithm for solving \eqref{eq:monotone_inv_mfg}.

\begin{algorithm}
\caption{Inverse method for the nonlocal mean-field game system
}\label{alg:short}
\begin{flushleft}
    \hspace*{\algorithmicindent} \textbf{Input:} $(\rho_{0,i}, g_i, \tilde{r}_{B,i})$ for $i = 1, ..., N$, $(\kappa_0, \mu_0)$\\
    \hspace*{\algorithmicindent} \textbf{Output:} $(\kappa^n, \mu^n)$ for $n = 1,...,\mathcal{N}_{\text{max}}$
\end{flushleft}
\begin{algorithmic}
    \While {iteration  $n<\mathcal{N}_{\text{maximal}}$}
    \State{1.Update for the adjoint problem:}
\State{   \hskip0.5em compute $( \lambda^{n+1}_{\rho_i}, \lambda^{n+1,}_{m_i} )$ use \eqref{alg:adjoint_1}\eqref{alg:adjoint_2} for $i = 1,...,N$.
    }
    \State{2. Update for the inverse problem:}
    \State{ \hskip0.5em compute $(\kappa^{n+1},\mu^{n+1})$ use \eqref{alg:inv_parameters}}
    \State{3. Update forward problem:}
    \State{ \hskip0.5em compute $(\rho_i^{n+1},m_i^{n+1},\phi_i^{n+1},a_i^{n+1})$ use \eqref{alg:update forward_mfg} for $i = 1,...,N$.}
    \State{$n \gets n+1$}
    \EndWhile 
\end{algorithmic}
\end{algorithm}

\subsection{Stabilizing techniques}
	
	Here, we discuss key numerical strategies for stabilizing Algorithm \ref{alg:short}.
	We refer to Appendix 2 (in particular, \textbf{Algorithm 2}) for more implementation details.
	
	While the change of $\kappa^{n+1}$ is made from the accumulation of all measurement events (through $( \lambda^{n+1}_{\rho_i}, \lambda^{n+1}_{m_i})$), there is sometimes unexpected change of $\kappa^{n+1}(x)$ that makes the algorithm highly unstable. For instance, there may be a large $\kappa^{n+1}(x)$ at a single grid point. Moreover, we are using harmonic expansion method to update $( \lambda^{n+1}_{\rho_i}, \lambda^{n+1}_{m_i})$, which causes large variances of $\kappa(x)$ along the boundary $\partial \Omega$. Therefore, we add a cut-off function and a convolution kernel to the step \eqref{alg:inv_parameters} to have a smoother change in $\kappa^{n+1}(x)$ in space. Specifically, we have
	\begin{equation*}
	\begin{cases}
	 \tilde{{\kappa}}^{n+1}  &=  \tilde{{\kappa}}^{n} +  \mathcal{T}_{mask}({\kappa}^{n+1, \text{temp}},\kappa_0)  - \kappa^{n}\\
	  \kappa^{n+1} &=  \max \left\{ \varepsilon_1, \tilde{{\kappa}}^{n+1}* \psi \right\}
	 		\end{cases},	    
	 		\end{equation*}
		where $\mathcal{T}_{mask}$ is a cut--off function that truncates the change of $\kappa$ near $\partial \Omega$ given by
		\begin{equation*}
		    \mathcal{T}_{mask}({\kappa},\kappa_0)(x) =\xi (x) (\kappa - \kappa_0)(x) + \kappa_0(x),
		\end{equation*}
		for a function $\xi(x)$ vanishing near $\partial \Omega$. As for the convolution 
			\begin{equation*}\kappa^{n+1}(x)  =  \max \left\{ \varepsilon_1, \int_{\Omega'}\tilde{{\kappa}}^{n+1}(y) \psi(x-y)dy \right\},	\end{equation*} where the convolution kernel $\psi(x)$ satisfies $\int_{\Omega} \psi(x) dx = 1$.

On the other hand, after the inverse problem parameters $(\kappa^{n+1},\mu^{n+1})$ are updated, we get	a new pair of parameters for a set of mean-field game problems. It is unclear whether starting from $(\rho_i^{n},m_i^{n},\phi_i^{n},a_i^{n})$ and taking the update rule \eqref{alg:update forward_mfg} once produces physical solutions for the new mean-field game system due to highly nonlinear dependence of the solution on the system parameters. Therefore, instead of preforming one iteration for the forward problem, we apply the PDHG algorithm for the forward problem until its error reaches a preset tolerance. More specifically, at every iteration $n$, with new system parameters $(\kappa^{n+1},\mu^{n+1})$, we use $(\rho_i^{n},m_i^{n},\phi_i^{n},a_i^{n})$ as an initial guess and calculate the mean-field game solution accurately so that the primal--dual gap is smaller than residual the preset tolerance.

	
	
	\section{Numerical examples} \label{sec:num_examples}

	This section demonstrates the efficiency and robustness of the inverse mean-field game algorithm with three examples. We also discuss details on the rule we used to choose the best reconstruction parameters.
	
	\subsection{Numerical implementation details}
	In this section, we present several numerical examples to illustrate the effectiveness
	of the new algorithm for the reconstruction of parameters in the mean-field game problem.
	
	We consider the spatial-time domain $\Omega'\times [0,T] = [-1,1]^2 \times [0,1] $. In the following examples, the partial boundary measurements are taken along the domain $\Omega = [-0.5,0.5]^2$, we refer as $\partial \Omega$. The Figure \ref{fig:example_domain} gives an example of the forward measurement event.

	In order to collect our observed data of the forward problem, we solve a set of mean-field game problem \eqref{eq:nonlocal_infsup} with given $(\rho_{0,i},g_i)$ and $(\kappa,\mu)$ by finite difference method with a mesh of size $(0.05,0.04)$ in space-time.
	Each problem is solved via primal-dual optimization approach with primal-dual gap $e_{tol}<2e-3$.
	The initial density function $\rho_{0,i}$ is the average of two Gaussian functions with centers ${x}_{G} \in \Omega' \backslash \Omega$.
	The final cost function $g(x)$ is smooth and has a smaller value around a single point $x_{g} \in \Omega' \backslash \Omega$ such that densities are concentrated in the neighborhood of $x_{g}$ at the final time.
	We want to point out that there is room to improve the initial density function and final cost function choices.
	We choose this set of $(\rho_{0,i},g_i)$ to ensure that the density's movement covers the domain $\Omega$ as completely as possible. 
	We also expect the nonlocal interaction among agents to be better reflected at the partial boundary measurements by setting the initial density as two Gaussians rather than one.
	
	We only take $16$ forward measurement events for each of the following numerical examples. 
	The partial boundary measurement means that we only collect the $\rho,m$ along the boundary $\partial \Omega$ in each event.
	Therefore, the resulting inverse problem is severely ill-posed.

	To test the robustness of our reconstruction algorithm, we add some random noise to the measurements as follows:
	\begin{equation}
	\left(  \rho, m \cdot n \right)^\delta(t_i,x_j) =  \left( (1+ \epsilon_n \delta_{ij,1}) \rho, (1+ \epsilon_n \delta_{ij,2}) m \cdot n \right) (t_i,x_j)\,, \label{noise}
	\end{equation}
	where $\{ (t_i,x_j) \}_{i=1,..,I,j=1,...,J} \in [0,T] \times \partial \Omega^+$ represents sampling points on the measurement boundary $[0,T] \times \partial \Omega^+$,
	$\{ \delta_{ij,1}, \delta_{ij,2} \}_{i=1,..,I,j=1,...,J}$ are i.i.d. random variables uniformly distributed on the interval $[-0.5,0.5]$ and $\epsilon_n$ corresponds to the noise level in the data, which is always set to be $\epsilon_n = 10 \%$ in all our examples.

 	From the noisy observed data $\{ \left(  \rho, m \cdot n \right)^\delta(t_i,x_j) \}_{i=1,..,I,j=1,...,J}$ on the sampling points of the measurement surface, we then use the algorithm to reconstruct the forward problem parameters $(\kappa,\mu)$. Recall that we paramatrized the running cost $L(x,v) := \frac{1}{2 \kappa(x)} |v|^2 $ by $\kappa$ and nonlocal kernel $K(x,y) :=\zeta(x;\mu,\omega)\cdot \zeta(y;\mu,\omega)$ by $\mu$. Since we aim at recovering the model on a given domain with fixed grid points, we fix the choice of $\omega$, and only seek sparse recovery of $\mu$.
	
	In the following examples, we use a set of parameters uniformly, without tuning.
	$\gamma_c = 0.2, \gamma_{\mu} = 0.1, \alpha_c =0.1,\alpha_{\mu} = 0.1, \alpha_{\lambda} = 10.$ We set the lower-bound projection parameter $\varepsilon_1=\kappa_c,$ this is based on the additional assumption of the model parameters that $\kappa(x) \geq \kappa_c\; \text{for}\; x \in \Omega'$. The projection parameter for kernel coefficient is $\varepsilon_{2} = 1$.
	
	To account for unknown ground truth of the model parameters, we introduce $\operatorname{Res}$ to quantify the quality of the reconstructed parameters. 
	\begin{align*}
	\operatorname{Res}^n = \sum_i \int_{\partial \Omega^+} \left(\|\rho_i^{n}  - \text{pr}_1 \tilde{r}_{B,i}\|^2 + \|  m_i^{n} \cdot n - \text{pr}_2 \tilde{r}_{B,i}  \|^2\right),
	\end{align*}
	where $\rho_i^n, m_i^{n}$ are the solution of the forward mean-field game problem with $i$-th choice of initial density and final cost function with the reconstructed parameter $(\kappa^n,\mu^n)$ at the $n$-th iteration of the algorithm. 
	The boundary residual $\operatorname{Res}$ measures how much the new boundary measurements of the mean-field game model with the recovered parameters deviate from the given partial measurements. If $(\kappa,\mu) =(\kappa_{true},\mu_{true})$, we would expect that $\operatorname{Res}$ is close to  $0$.
	Therefore, we pick the reconstructed parameters at $n_{opt}$-th iteration by taking
	\begin{align*}
	n_{opt} =\argmin\limits_n \operatorname{Res}^n,\\
	(\kappa_{opt},\mu_{opt}) = 	(\kappa^{n_{opt}},\mu^{n_{opt}}). 
	\end{align*} 
	When we implemented the algorithm, we observed that the quantity $\operatorname{Res}^n$ first decreased then increased with respect to the iteration. We also observed that with large enough number of iterations, (for example, $1500$), the inverse problem is contaminated and the reconstruction of mean-field game coefficients are very bad. In the following examples, we take fixed number of iterations $N =1500$ for the inverse algorithm, and pick the reconstructed model parameters accordingly.
	
	\subsection{Example 1} This example tests a running cost $\kappa(x)$ with a bump at point $(0.25,0.25)$, which means the density that travel crossing near this point has a lower cost than other routes.
	The density are also expected to accelerate when they travel across this point.
	The nonlocal kernel $K(x,y)$ is constructed via a Gaussian function plus some sparse terms in forms of $\mu_k^2\cos(\omega_k \cdot x)$.
	The nonlocal kernel, in general, penalizes being too concentrated.
	The amplify of certain Fourier frequencies determines the agents' particular interaction preferences.
	Specifically, we have the following: 
	\begin{align*}
	\kappa(x) & = 2 + 4 \exp{\left(-\frac{(x_1-0.25)^2 + (x_2-0.25)^2}{0.1^2}\right)},\\
		\kappa_0(x) & = 2,\\
	K(x,y) & = K_0(x,y) + K_s(x,y) + k_0,\\
	K_0(x,y) & = \frac{1}{5} \exp{\left(-\frac{1}{2}\frac{x^2 + y^2}{0.4^2}\right)}, \\
	K_s(x,y) & =  0.2094^2 \left(\cos{\left(\pi (x_1-y_1)\right)} + \cos{\left(\pi (x_2 -y_2)\right)}\right)\\ & + 0.2613^2\left(\cos{\left(\pi (x_1-y_1)+ \pi (x_2 -y_2)\right)} + \cos{\left(-\pi (x_1-y_1)+ \pi (x_2 -y_2)\right)}\right),
	\end{align*}
	where $x = [x_1,x_2], y = [y_1,y_2]$.
	We have $\mu_0$, which represents $K_0(x,y)$ via the expansions form  $\mu_k^2\cos(\omega_k \cdot x)$,  known. The variable $k_0$ is a given constant value that makes the kernel integration $\int\int K(x,y) dx dy = 1$. Varying this constant corresponding to changing the coefficient of the zero Fourier mode $(0,0)$. This constant $k_0$ does not change the intensity of repulsion effect among the agents, since $\int k_0 \rho(y) dy = k_0\int \rho_0(y) dy$ is uniform over the domain $\Omega'$.
	With $K(x,y)=\zeta(x;\mu,\omega)\cdot \zeta(y;\mu,\omega)$ for $\mu,\omega \in C_{odd=even}$, we omit the even entries (eg. $\mu_{k,2}, \omega_{k,2}$) and express the kernel $K_s$ as follows:
	\begin{align*}
	\mu_s &= (0.2094,0.2094,0.2613,0.2613), \\
	\omega_s & = ((\pi,0),(0,\pi),(\pi,\pi),(-\pi,\pi)).
	\end{align*}
Here, we also assume that $\kappa_c =2$ and $\kappa(x) = 2$ for $ x \in  \Omega' \backslash \Omega$ is known.
	\begin{figure}[htbp!]
		\includegraphics[width=0.950\textwidth,trim=0 0 0 00, clip=true]{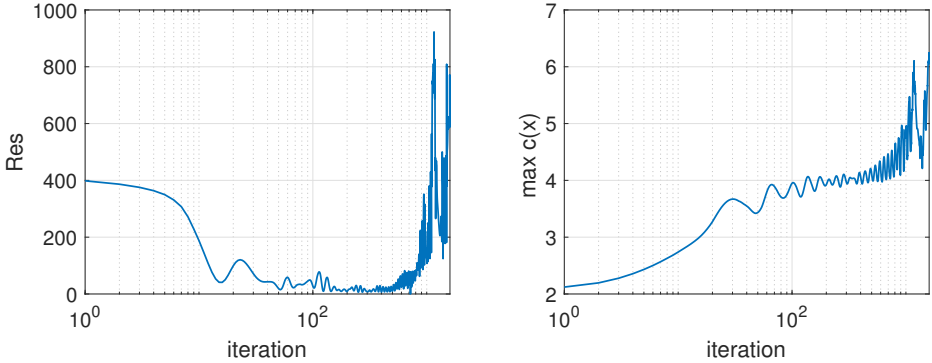}
		\caption{The residual $\operatorname{Res}^n$ and the $\max_x \kappa^n(x)$ at $n$-th iteration.}
		\label{fig:eg1_res}
	\end{figure}
	\begin{figure}[htbp!]
		\includegraphics[width=0.950\textwidth,trim=20 20 20 25, clip=true]{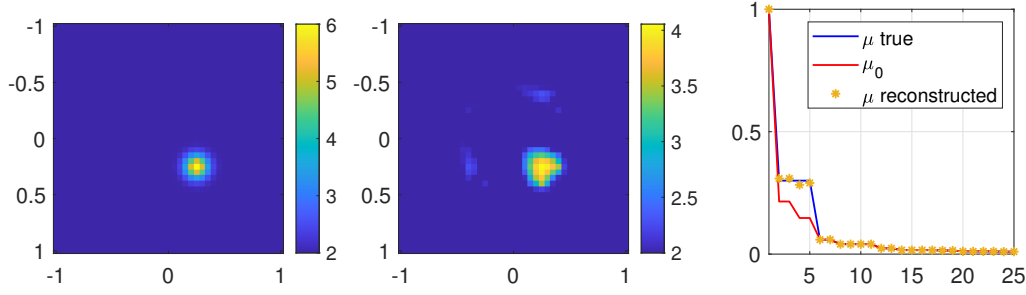}
		\caption{From left to right: the true running cost $\kappa(x)$; the reconstructed running cost $\kappa_{opt}(x)$ at iteration $n_{opt}$; the coefficient representation of nonlocal kernel $K(x,y)$ in vector form, where $x$-axis represents different Fourier mode $\omega$ and the $y$-axis corresponds to the coefficients $\mu$.  }
		\label{fig:eg1}	
	\end{figure}
	
	Given $(\kappa_0,\mu_0)$ and the noisy partial boundary measurements with corresponding event parameters $(\rho_{0,i}, g_i, \tilde{r}_{B,i})$, we apply our inverse algorithm.The results are shown in Figure \ref{fig:eg1_res},\ref{fig:eg1}.
	In Figure \ref{fig:eg1_res}, we plot the residual $\operatorname{Res}^n$ and the $\max_x \kappa^n(x)$ along the iteration. We see that the residual oscillates and decreases first, then bounces back and increases.
	In Figure \ref{fig:eg1}, we show the reconstruction of model parameters by taking $n_{opt} =\argmin_n \operatorname{Res}^n,
	(\kappa_{opt},\mu_{opt}) = 	(\kappa^{n_{opt}},\mu^{n_{opt}}).$ We see that the reconstructed $\kappa_{opt}(x)$ has a single bump sits near $(0.25,0.25)$. The shape of the bump is not as sharp as the ground truth $\kappa$. The maximal value of running cost $\max_x \kappa_{true}(x) = 6$; while $\max_x \kappa_{opt} (x) = 4$.
	As for the non-local kernel, we have $\mu_{opt}$ nicely reconstructed, where $\mu_{true} = \mu_0 + \mu_s \approx \mu_{opt}$.  
	This example shows that our inverse algorithm is robust to noise and can recover the model parameters $(\kappa,\mu)$ simultaneously.
\subsection{Example 2}
	In this example, we make the $\kappa(x)$ more complicated by having two bumps sitting diagonally.
	We except that if the density travels across these two bumps, it will accelerate twice. The model set-up is as follows:
	\begin{align*}
	\kappa(x) & = 2 + 4 \exp{\left(-\frac{(x_1+0.25)^2 + (x_2-0.25)^2}{0.1^2}\right)} + 4 \exp{\left(-\frac{(x_1-0.25)^2 + (x_2+0.25)^2}{0.1^2}\right)},\\
	\kappa_0(x) & = 2, \kappa_c =2,\\
	K(x,y) & = K_0(x,y) + K_s(x,y)+ k_0,\\
	K_0(x,y) & = \frac{1}{5} \exp{\left(-\frac{1}{2}\frac{x^2 + y^2}{0.4^2}\right)} ,\\
	\mu_s &= (0.3374, 0.3374, 0.2942, 0.2942) ,\\
	\omega_s & = ((\pi,0),(0,\pi),(2 \pi, 0),( 0, 2\pi)).
	\end{align*}
	\begin{figure}[htbp!]
		\includegraphics[width=0.950\textwidth,trim=20 20 20 40, clip=true]{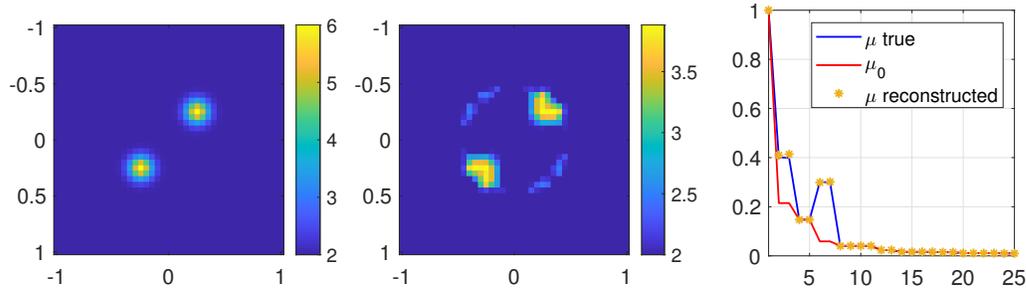}
		\caption{From left to right: the true running cost $\kappa(x)$; the reconstructed running cost $\kappa_{opt}(x)$ at iteration $n_{opt}$; the coefficient representation of nonlocal kernel $K(x,y)$ in vector form. }
		\label{fig:eg2}	
	\end{figure}

	We can see from the Figure \ref{fig:eg2} that recovered bumps are well separated, and their locations are accurately captured. Reconstructed bumps are more spread compared to the ground truth, and there is some noise on upper left and bottom right corners of the domain $\Omega'$. The nonlocal kernel is reconstructed nicely as shown in Figure \ref{fig:eg2}(right). A precise sparse representation of $K_s(x,y)$ is recovered. Considering the severe ill-posedness of the inverse problem with $10\%$ multiplicative noise added to the boundary measurements, the reconstruction quality is quite satisfactory. 
	
	\subsection{Example 3}
	In this example, we modify the $\kappa(x)$ by having two bumps sitting in parallel.
	Similar to the Example $2$, the density would prefer to move crossing these bumps. We set the nonlocal kernel with $K_s$ containing Fourier modes with higher frequency. 
	\begin{align*}
	\kappa(x) & = 2 + 4 \exp{\left(-\frac{(x_1-0.25)^2 + (x_2-0.25)^2}{0.1^2}\right)} + 4 \exp{\left(-\frac{(x_1-0.25)^2 + (x_2+0.25)^2}{0.1^2}\right)}.\\
		\kappa_0(x) & = 2, \kappa_c =2,\\
	K(x,y) & = K_0(x,y) + K_s(x,y)+ k_0,\\
	K_0(x,y) & = \frac{1}{5} \exp{\left(-\frac{1}{2}\frac{x^2 + y^2}{0.4^2}\right)} ,\\
	\mu_s &= ( 0.2973,  0.2973,  0.2973,  0.2973) ,\\
	\omega_s & = ((2\pi, -\pi),(2\pi, \pi),(\pi, 2\pi),( \pi, -2\pi)).
	\end{align*}
	\begin{figure}[htbp!]
		\includegraphics[width=0.950\textwidth,trim=20 20 20 40, clip=true]{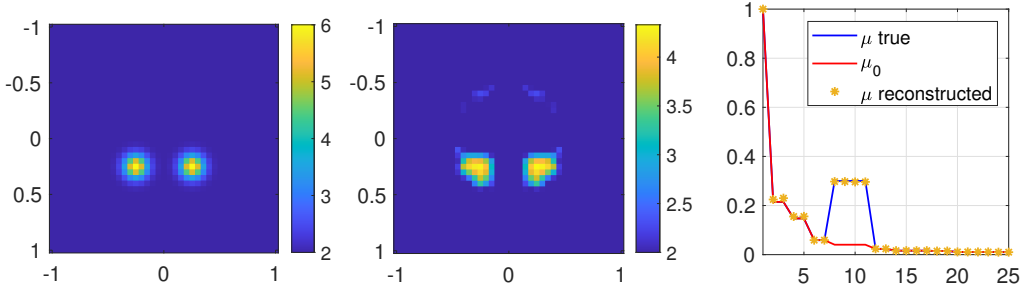}
		\caption{From left to right: the ground true running cost $\kappa(x)$; the reconstructed running cost $\kappa_{opt}(x)$ at iteration $n_{opt}$; the coefficient representation of nonlocal kernel $K(x,y)$ in vector form. }
		\label{fig:eg3}	
	\end{figure}
	
	We have the reconstruction result shown in Figure \ref{fig:eg3}. The two parallel sitting bumps are well separated and located with reasonable accuracy. Again, the bumps are diffused with some noise near the upper boundary of $\Omega$. The nonlocal kernel is recovered very nicely.  
	
\section{Conclusion}\label{sec:conclusion}
	
	In this paper, we formulate a new class of inverse mean-field game problems given only partial boundary measurements. A novel model recovery algorithm is proposed based on the saddle point formulation of MFGs. We demonstrate the robustness and effectiveness of the numerical inverse algorithm with several examples, where the MFG model parameters are reconstructed accurately. Our algorithm can be further generalized to other inverse problems with saddle point structure in the forward problem. 

\bibliographystyle{plain}
\bibliography{literature}

\begin{appendices}

\section{The KKT conditions of the inverse mean-field game problem}

We denote by
\begin{align*}
\mathcal{L}(\cdot) : = \mathcal{L} \bigg(  (\rho_{0,i}, g_i ),  ( \rho_i, m_i), (a_i,\phi_i)  , (\kappa, \mu)  \bigg),
\end{align*}
the function $\mathcal{L}$ defined in Equation \ref{eq:nonlocal_infsup} for simplicity. The KKT conditions for the inverse mean-field game problem are then
	\begin{eqnarray*}
	\begin{cases}
		\Pi_{B,(\rho,m)}^* [ \Pi_{B,(\rho,m)} (\rho_i, m_i) - \tilde{r}_{B,i}  ] + \partial^2_{(\rho_i,m_i),(\rho_i,m_i)}
		\mathcal{L}(\cdot)  \lambda_{(\rho_i,m_i)}   = 0, \quad i = 1, ...,N,\\
		\Pi_{B,(a,\phi)}^* [ \Pi_{B,(a,\phi)} (a_i,\phi_i) - \tilde{s}_{B,i}  ] - \partial^2_{(a_i,\phi_i),(a_i,\phi_i)}  \mathcal{L} (\cdot)   \lambda_{(a_i,\phi_i)}   = 0, \quad i = 1, ...,N,\\
		\partial_{(\kappa, \mu)} R_1(\kappa, \mu) +  \partial_{(\kappa, \mu)} R_2 (\kappa, \mu) 
		+ \sum_{i=1}^N \left \langle  \partial_{(\kappa, \mu)} \partial_{(\rho_i,m_i)}
		\mathcal{L} (\cdot)  , \lambda_{(\rho_i,m_i)} \right \rangle,    \\
		- \sum_{i=1}^N  \left \langle \partial_{(\kappa, \mu)}  \partial_{(a_i,\phi_i)}  \mathcal{L}(\cdot)  , \lambda_{(a_i,\phi_i)} \right \rangle    = 0, \\
		\partial_{(\rho_i,m_i)}  \mathcal{L} (\cdot)   = 0, \quad i = 1, ...,N,\\
		-  \partial_{(a_i,\phi_i)}  \mathcal{L}(\cdot)   = 0, \quad i = 1, ...,N.
	\end{cases}
\end{eqnarray*}
Furthermore, the derivatives of $\mathcal{L}$ are given by
{
		\begin{eqnarray*}
			\partial_{(\rho_i,m_i)} \mathcal{L}  &=& \left( \partial_t \phi+ \nu \Delta \phi - \frac{1}{2 \kappa} \frac{|m|^2}{\rho^2} + a \cdot \zeta(\cdot;\mu,\omega) ,   \nabla \phi  + \frac{1}{\kappa } \frac{m}{\rho}  \right)\\
			\partial_{(a_i,\phi_i)}  \mathcal{L}  &=& \left( - a + \int_{\Omega'} \rho(y) \,  \zeta( y; \mu, \omega)dx ,  - \partial_t \rho + \nu \Delta \rho - \nabla \cdot m \right) \\
			\partial^2_{(\rho_i,m_i),(\rho_i,m_i)} \mathcal{L}  &=&  \begin{pmatrix}  \frac{1}{\kappa} \frac{|m|^2}{\rho^3} &  - \frac{1}{\kappa} \frac{m}{\rho^2} \\  - \frac{1}{\kappa} \frac{m}{\rho^2} &    \frac{1}{\kappa \rho} I  \end{pmatrix}\\
			\partial^2_{(a_i,\phi_i),(a_i,\phi_i)}  \mathcal{L}  &=& \begin{pmatrix} I & 0 \\  0  & 0 \end{pmatrix} \\
			\partial_{(\kappa,\mu)} \partial_{(\rho_i,m_i)}
			\mathcal{L} 
			&=& \begin{pmatrix}  \frac{1}{2 \kappa^2(x)} \frac{|m|^2}{\rho^2} & - \frac{1}{\kappa^2(x)} \frac{m}{\rho}  \\  
				\Lambda_1 (\omega,\cdot) \, a
				& 0 \end{pmatrix} \\
			\partial_{(\kappa,\mu)}  \partial_{(a_i,\phi_i)}  \mathcal{L} 
			&=& \begin{pmatrix} 0 & 0 \\   
				\int_{\Omega'}\rho(y) \, \Lambda_1 (\omega,y) dy
				& 0\end{pmatrix} 
		\end{eqnarray*}
	} 
where the variable $\Lambda_1 (\omega,x)$ is defined in Section \ref{sec:alg_details:4-1-2}.

\section{The inverse algorithm}

\begin{algorithm}
\caption{A detailed inversion algorithm for the nonlocal mean-field game system
}\label{alg:long}
\begin{flushleft}
    \hspace*{\algorithmicindent} \textbf{Input:} $(\rho_{0,i}, g_i, \tilde{r}_{B,i})$ for $i = 1, ..., N$, $(\kappa_0, \mu_0)$\\
    \hspace*{\algorithmicindent} \textbf{Output:} $(\kappa^n, \mu^n)$ for $n = 1,...,\mathcal{N}_{\text{max}}$
\end{flushleft}
\begin{algorithmic}
    \While {iteration  $n<\mathcal{N}_{\text{maximal}}$}
    \State{1.Update for the adjoint problem by computing $( \lambda^{n+1}_{\rho_i}, \lambda^{n+1,}_{m_i} )$}
\State{   \hskip0.5em \(
\begin{cases}
\begin{cases}
\Delta  \eta_i^{n+1} = 0 \text{ in } \Omega\\
 \eta_i^{n+1}  =  \rho_i^n  - \text{pr}_1 \tilde{r}_{B,i}  \text{ on } \partial \Omega^+
\end{cases} \\
\begin{cases}
\Delta  \eta_i^{n+1} = 0 \text{ in } \Omega' \backslash \Omega\\
 \eta_i^{n+1}  =  \rho_i^n  - \text{pr}_1 \tilde{r}_{B,i}  \text{ on } \partial \Omega^- \\
 \eta_i^{n+1}  = 0  \text{ on } \partial \Omega'
\end{cases}\\
\begin{cases}
\Delta \xi_i^{n+1} = 0 \text{ in } \Omega\\
\partial_n  \xi_i^{n+1}   =   m_i^n \cdot n - \text{pr}_2 \tilde{r}_{B,i}  \text{ on } \partial \Omega^+
\end{cases} \\
\begin{cases}
\Delta \xi_i^{n+1} = 0 \text{ in } \Omega' \backslash \Omega\\
\partial_n  \xi_i^{n+1}   =   m_i^n \cdot n + \text{pr}_2 \tilde{r}_{B,i}  \text{ on } \partial \Omega^- \\
\partial_n  \xi_i^{n+1}   =  0 \text{ on } \partial \Omega'
\end{cases} 
\end{cases}
\) \\ }
\State{   \hskip0.5em \(
\begin{cases}
&\begin{pmatrix} \lambda^{n+1,\text{temp}}_{\rho_i}  \\ \lambda^{n+1,\text{temp}}_{m_i} \end{pmatrix} =\\ 
& \left [ I + \alpha^n_{\lambda_{(\rho_i,m_i)}} \rho_i^n \begin{pmatrix}  \frac{1}{\kappa^n} \frac{|m^n_i|^2}{(\rho^n_i)^3}   - \frac{1}{\kappa^n} \frac{m^n_i}{(\rho^n_i)^2} \\  - \frac{1}{\kappa^n} \frac{m_i^n}{(\rho^n_i)^2}     \frac{1}{\kappa^n \rho_i^n} I  \end{pmatrix} \right]^{-1}
\left( \begin{pmatrix} \lambda^{n}_{\rho_i}  \\ \lambda_{m_i^{n}} \end{pmatrix} - \alpha^n_{\lambda_{(\rho_i,m_i)}} \rho_i^n  \begin{pmatrix} \eta^{n+1}_i  \\ \nabla \xi^{n+1}_i \end{pmatrix} \right)\\
&\begin{pmatrix} \lambda^{n+1}_{\rho_i}  \\ \lambda^{n+1}_{m_i} \end{pmatrix}
= 2
\begin{pmatrix} \lambda^{n+1,\text{temp}}_{\rho_i}  \\
	\lambda^{n+1,\text{temp}}_{m_i} \end{pmatrix}
-\begin{pmatrix} \lambda^{n, \text{temp}}_{\rho_i}  \\ 
\lambda^{n, \text{temp}}_{m_i} \end{pmatrix} 
\end{cases}
\)   }
    \State{2. Update for the inverse problem by computing $(\kappa^{n+1},\mu^{n+1})$}
    \State{ \hskip0.5em \(
\begin{cases}
\begin{cases}
{\kappa}^{n+1, \text{temp}}  &= S_{\alpha_{\kappa}^n  \gamma}\bigg(  2 \kappa^{n} - \tilde{{\kappa}}^{n} - \alpha_{\kappa}^n  
\sum_{i=1}^N   \Lambda_{\kappa}(\kappa^n,m_i^n,\rho_i^n,\lambda_{\rho_i}^{n+1},\lambda_{m_i}^{n+1} )    - \kappa_0  \bigg) + \kappa_0 \\
	 \tilde{{\kappa}}^{n+1}  &=  \tilde{{\kappa}}^{n} +  \mathcal{T}_{mask}({\kappa}^{n+1, \text{temp}},\kappa_0)  - \kappa^{n} \\
	 	  	  \kappa^{n+1} &=  \max \left\{ \varepsilon_1, \tilde{{\kappa}}^{n+1}* \psi \right\}
\end{cases} \\
\begin{cases}
{\mu}^{n+1,\text{temp}}   &=
 S_{\alpha_{\mu}^n  \gamma}\bigg(  2\Pi^*_{\mu}(\mu^n) - \tilde{\mu}^n- \alpha_{\mu}^n  
\sum_{i=1}^N  \Lambda_{\mu}(\lambda_{\rho_i}^{n+1},a_i^n)  
 - \Pi^*_{\mu}(\mu_0)  \bigg) + \Pi^*_{\mu}(\mu_0)  \\
\tilde{\mu}^{n+1}  & = \tilde{\mu}^{n} + {\mu}^{n+1,\text{temp}} - \mu^{n} \\
\mu^{n+1} & = \Pi_{\mu} (\min \{ \varepsilon_2, \tilde{\mu}^{n+1}\} ) 
 
 \end{cases} 
\end{cases} \)\\}
    \State{3. Update forward problem by computing $(\rho_i^{n+1},m_i^{n+1},\phi_i^{n+1},a_i^{n+1})$ with the forward mean-field game algorithm for $i = 1,...,N$.}
    \State{ \hskip0.5em Apply the iterative Algorithm \ref{alg:forward} given input $(\rho_{0,i}, g_i)$, $(\kappa^{n+1}, \mu^{n+1})$, $e_{tol}$ with initial guess $(\rho^n_i,m^n_i,\phi^n_i,a^n_i)$, and assign
    \begin{equation*}
    (\rho_i^{n+1},m_i^{n+1},\phi_i^{n+1},a_i^{n+1}) := (\rho_i^{*},m_i^{*},\phi_i^{*},a_i^{*})
    \end{equation*}}
    \State{$n \gets n+1$}
    \EndWhile 
\end{algorithmic}
\end{algorithm}

\newpage

\section{The forward algorithm for nonlocal mean-field games }

\begin{algorithm}
\caption{Iterative algorithm for the nonlocal mean-field game system
}\label{alg:forward}
\begin{flushleft}
    \hspace*{\algorithmicindent} \textbf{Input:} $(\rho_{0}, g)$, $(\kappa, \mu)$, a set of initial guess $(\rho^0,m^0,\phi^0,a^0)$, $e_{tol}$, a set of stepsizes $(\alpha_{\rho}^j,\alpha_{m}^j,\alpha_{\phi}^j,\alpha_{a}^j)$\\
    \hspace*{\algorithmicindent} \textbf{Output:} $(\rho^*,m^*,\phi^*,a^*)$
\end{flushleft}
\begin{algorithmic}
    \While {iteration  $j<\mathcal{J}_{\text{max}}$ and primal-dual gap $\text{PD}(\rho^j,m^j, \phi^{j} , a^{j})\geq e_{tol}$}
    \State{ 
    \( \begin{cases}
( \rho^{j+1}  , m^{j+1} ) &= 
\text{argmin}_{(\rho,m)} \bigg \{ \mathcal{L} \bigg ( (\rho_{0}, g ) , ( \rho, m), (a^j,\phi^j)  , (\kappa,\mu)  \bigg) \\
 &+ \frac{1}{2 \alpha_{\rho}^j} \|\rho^{j} - \rho \|_{L^2_x,t}^2  + \frac{1}{2 \alpha_{m}^j} \|m^{j} - m \|_{L^2_{x,t}}^2  \bigg \}
\\
( \phi^{j+1 {, \text{temp}} } , a^{j+1 {, \text{temp}} } ) &= \text{argmin}_{(a,\phi)}
 \bigg\{ - \mathcal{L} \bigg ( (\rho_{0}, g ), ( \rho^{n { +1 }}, m^{n { +1 } }), (a,\phi) ,(\kappa,\mu)  \bigg)  \\
 &+ \frac{1}{2 \alpha_{\phi}^j } \| \phi^{j} -  \phi \|_{H^1_{x,t}}^2  + \frac{1}{2 \alpha_{a}^j } \| a^{j} -  a \|_{L^2_{t}}^2
 \bigg\} \\
  ( \phi^{j+1} , a^{j+1} )  & = 2 ( \phi^{j+1,  \text{temp} } , a^{j+1, \text{temp} } )  - ( \phi^{j} , a^{j} ) 
\end{cases} \)}
    \State{$j \gets j+1$}
    \EndWhile 
    \State \Return {$(\rho^j,m^j, \phi^{j} , a^{j})$}
\end{algorithmic}
\end{algorithm}
\end{appendices}

\end{document}